\documentstyle[11pt]{article}

\newtheorem{th}{Theorem}[section]
\newtheorem{de}[th]{Definition}

\newtheorem{lem}[th]{Lemma}

\newtheorem{re}[th]{Remark}

{

\begin{document}

\large
\title {\bf ON SOME ESTIMATES FOR PROJECTION OPERATOR IN BANACH SPACE}
\author{\bf Ya. I. Alber$^{*}$ and A. I. Notik \\
Department of Mathematics\\
Technioin-Israel Institute of Technology\\
Haifa 32000, Israel}
\date{}
\maketitle

{\bf Abstract} -
Previously unknown estimates of  uniform  continuity  of  projection
operators in Banach space have been obtained.  They can be  used  in  the
investigations of  approximation  methods,  in  particular,
the  method  of quasisolutions,  methods  of  regularization  and  penalty
functions,  for  solving
nonlinear problems on exact and peturbed sets (see \cite{a1,a2}).\\

\section{Introduction}
\setcounter{equation}{0}

Metric projection operators on convex closed sets in Hilbert
and Banach spaces are widely used in
functional and numerical analysis, theory of optimization
and theory of approximation.  Terefore their properties
are actively studying .

Metric projection operators  can be defined in a same way
in Hilbert and Banach spaces.

Let $B$ be a real, uniformly convex and uniformly  smooth  (reflexive)
Banach space with $B^*$ its dual space \cite{lt} and $J: B \rightarrow B^* $ a
normalized  duality mapping,
$\Omega$  a closed convex set in $B,\; <w,v>$ a dual product
in $B$, i.e. pairing  between
$w \in  B^*$ and $v \in  B \;((\cdot,\cdot)$ is inner product in Hilbert
space);
the sigs $||\cdot||, \;||\cdot||_{B^{*}}$ and $ ||\cdot||_{H} $ denote
the norms in Banach space $B$, Banach space $B^{*}$ and Hilbert space $H$,
respectively.
Let $P_{\Omega }x$ and $P_{\Omega }y$ be projections of
elements $x$ and $y$ onto $\Omega $ in the sense of best approximation:
\begin {de}  \label{gpo}
The operator $P_{\Omega}: B \rightarrow {\Omega} \subset B$
is called to be metric projection operator in Banach space $B$ if it gives
the correspondence between an arbitrary point  $x \in  B$ and nearest
point
 $\bar x \in \Omega $ according to minimization problem
\begin{equation} \label{f1}
P_{\Omega}x = \bar x; \;\;\; \bar  x: ||x - \bar x|| =
\inf_{\xi \in \Omega} ||x - \xi||.
\end{equation}
\end {de}
------------------------\\
$^{*}$This research was supported in part by the Ministry of Science Grant
3481-1-91 and by the Ministry of Absorption Center for Absorption in Science.

It is well known that the metric projection operator in Hilbert space $ H$
is:\\
1) monotone (accretive), i.e. for all $x,y \in H $
$$(\bar x - \bar y, x-y) \ge 0, $$
2) nonexpansive, i.e. it satisfies the Lipschitz condition (with
the constant 1)
 $${ ||\bar x - \bar y||_H} \leq {||x- y||_H},
{\qquad }
\forall x, y \in H .$$
3) It gives an absolutely best
approximation for arbitrary elements from Hilbert space by the
elements of convex closed set $\Omega$ , i.e.
$${||\bar x  - \xi||}^2 _H \leq {||x-\xi||}^2 _H - {||x - \bar x||}^2_H
,
{\qquad }
\forall \xi \in\Omega ,$$
4) The operator $P_{\Omega}$ satisfies the following necessary
and sufficient condition
$$ (x - \bar x, \bar x -\xi) \ge 0,  {\qquad }        \forall
\xi \in\Omega  . $$
Mainly, these properties lead to a variety of applications of
the metric projection operator in Hilbert spaces.

Metric projection operators in Banach spaces do not have the properties 1)-3)
mentioned above . However, the property 4) for Banach spaces has been
established
in \cite{jl} (see \cite{gr} also) in a form
\begin{equation} \label{f0}
< J(x - \bar x), \bar x -\xi > \ge 0,  {\qquad }
\forall \xi \in\Omega.
\end{equation}
This property is often used in applications. A question
connected with  smoothness of the metric projection operators
in Banach spaces is not of less importance either.
It is known that in uniformly  smooth  Banach space
the operator $P_{\Omega}$ is always continuous but not  always uniformly
continuous. In the present paper, we prove the uniform continuty of
$P_{\Omega}$
with respect to changes of the argument $x$ (Theorem \ref{th1})
in uniformly  smooth and  uniformly  convex Banach spaces and
the uniform continuty of $P_{\Omega}$ with respect to changes
of set $\Omega$ (Theorem \ref{th2}) in uniformly convex Banach spaces. This
properties have been applied for investigation of stability of approximation
methods for solving equations and variational inequalities in Banach spaces
\cite{a1,a2}.

\section{Auxiliary lemmas}
\setcounter{equation}{0}
In the prooving  main results, we use some properties of duality mappings
which  were in most of known monographs devoted
to the theory of monotone and accretive operators. \\
i) The normalized  duality mapping $J$ is  monotone operator, i.e.
$$<Jx - Jy, x - y> \ge 0, {\qquad }  \forall \xi \in\Omega.$$
ii) In uniformly  convex Banach spaces, mapping $J$ is
uniformly monotone operator on each bouded set, i.e. for every $R>0$ there
exists a function $\psi_{R}(t) > 0, \psi_{R}(0) = 0$, such that for all
$x$ and $y$ in $B$, $||x|| \leq R$ and $||y|| \leq R,$
$$<Jx - Jy, x - y> \ge \psi_{R}(||x - y||) .$$
iii) In uniformly  smooth Banach spaces, mapping $J$ is
uniformly continuous operator on each bouded set, i.e. for every $R>0$ there
exists a continuous function $\omega_{R}(t) > 0, \;\omega_{R}(0) = 0$, such
that for all  $||x|| \leq R$ and $||y|| \leq R,$
$$||Jx - Jy||_{B^{*}} \leq \omega_{R}(||x - y||).$$
However,
during a big period of time  the functions $\psi_{R}(t)$ and $\omega_{R}(t)$
had no quantitative description. And only in 1984, such estimates have
been obtained in \cite{an1} (see also \cite{an2}). But in  \cite{an1,an2} we
did not proved their detailed proofs. In this Section we are intended to
fill this gap.
\begin {lem} \label{l1}
(cf.\cite{an1,an2}). Let $B$ be an uniformly  convex Banach spaces.
If $\delta _{B}(\epsilon )$ is modulus of the  convexity  of  space $B$
\cite{lt}
then the following estimate
$$<Jx - Jy, x - y> \ge (2L)^{-1}\delta _{B} (||x -  y||/C_{1}) ,$$
$$ C_1 = 2 \max \lbrace  1, \sqrt{(||x||^2 + ||y||^2)/2}\rbrace $$
is valid. Here $L $  is the constant from the Figiel's
iniquality (\ref{a14}),$\; 1 < L < 3,18 .$
\end {lem}
{\bf Proof}. In \cite{lt} it was shown that if $x \in B, y \in B $ and
$||x||^2 + ||y||^2 = 2$, then
$$||(x + y)/2||^2 \leq 1 - \delta_B (||x-y|| / 2)$$
where $\delta_B (\epsilon)$
is modulus of convexity of space $B$. It is known that in uniformly
convex Banach space, $\delta_B(\epsilon)$ is convex strictly
increasing function,
$0 \leq \delta_B (\epsilon) < 1,  \delta_B (0) = 0 .$
We denote $R^2_1 = 2^{-1}
(||x||^2 + ||y||^2)$ ($x$ and $y$ are not zero at the same time) and
introduce
new variables $ \tilde x = x/R_1$ and $ \tilde y = y/R_1$. Then
$$||\tilde x||^2 + ||\tilde y||^2 = R^{-2}_1(||x||^2 + ||y||^2) = 2 .$$
Therefor, for $\tilde x$ and $\tilde y$ the inequality
$$||(\tilde x + \tilde y)/2||^2 \leq 1 -
\delta_B (||\tilde x - \tilde y||/2)$$
is valid.

Let us return to the old variables $x$ and $y.$ We obtain
$$||(x + y)/2R_1||^2 \leq 1 - \delta_B (||x-y|| / 2R_1)$$
or
\begin{equation} \label{a3}
||(x + y)/2||^2 \leq (||x||^2 + ||y||^2)/2 - R^2_1 \delta_B (||x-y|| / 2R_1),
\end{equation}
$$R_1 = \sqrt{(||x||^2 +||y||^2)/2} .$$
Consider two cases:

1. $R_1 \ge 1.$ Then
\begin{equation} \label{a4}
||(x + y)/2||^2 \leq (||x||^2 + ||y||^2)/2 - \delta_B (||x-y|| / 2R_1).
\end{equation}

2. $R_1 \leq 1.$ From (\ref{a3}) and the inequality  (see \cite{f})
\begin{equation} \label{a14}
\epsilon ^2 \delta_B (\eta) \ge (4L)^{-1} \eta ^2 \delta_B (\epsilon),
\;\;\;\forall \eta \ge \epsilon >0
\end{equation}
we have
$$\delta_B (||x-y|| / 2R_1) \ge R^{-2}_1 (4L)^{-1} \delta_B (||x-y|| / 2) $$
and
\begin{equation} \label{a5}
||(x + y)/2||^2 \leq (||x||^2 + ||y||^2)/2 - (4L)^{-1} \delta_B (||x-y|| / 2).
\end{equation}
In vertue of $\min \lbrace 1, (4L)^{-1} \rbrace = (4L)^{-1},$
(\ref{a4}) and (\ref{a5}) if joined together give
$$||(x + y)/2||^2 \leq (||x||^2 + ||y||^2)/2 - (4L)^{-1}
\delta_B (||x-y|| / C_1),\;\; C_1 = 2\max \lbrace 1,R_1 \rbrace .$$
If $||x|| \leq R$ and $||y|| \leq R$, then $C_1 = 2\max \lbrace 2,R \rbrace$
is the absolute constant.

Let $\varphi (x) = ||x||^2/2.$ Then
\begin{equation} \label{a6}
\varphi ((x +y)/2) \leq 2^{-1} \varphi (x)  + 2^{-1} \varphi (y) -
(8L)^{-1} \delta_B (||x-y|| / C_1).
\end{equation}
In \cite{vnc} it was proven that
if a convex functional ${\phi(x)}$ difined on convex closed set $ \Omega $
satisfies the inequality
$${\phi({\frac 12}x + {\frac 12} y)} \leq {\frac 12} \phi (x) +
{\frac 12} \phi (y) -  \kappa (||x - y||),$$
where $\kappa (t) \ge 0,  \kappa (t_0) > 0 $ for some $t_0 > 0$,
then ${\phi(x)}$ is uniformly convex  functional with modulus of convexity
$\delta (t) = 2 \kappa (t)$, and
$${\phi(x)} \ge  \phi (y) + <l(y),x - y> + 2 \kappa (||x - y||) ,$$
\begin {equation} \label{f109}
 <l(x) - l(y),x - y> \ge 4 \kappa (||x - y||)
\end{equation}
for all $l(x) \in \partial \phi(x)$. Here $ \partial \phi(x)$ is the
set of all support functionals (the set of all subgradients) of $\phi(x)$ at
the point $x \in \Omega.$ \\
The formula  (\ref {a6}) then shows that the functional
$\varphi (x)$ is uniformly convex
with modulus of convexity
$$\delta (||x-y||) = (4L)^{-1} \delta_B (||x-y|| / C_1). $$
Let us apply  now (\ref {f109}). Because
$\varphi (x)$ is differentiable functional and
$\varphi'(x) = grad (||x||^2/2) = Jx$, then
\begin{equation} \label{a15}
<Jx - Jy, x - y> \ge (2L)^{-1} \delta_B (||x-y||/C_1).
\end{equation}
Lemma is proved.
\begin {lem} \label{l2}
(cf.\cite{an1,an2}). Let ${B^{*}}$ be an uniformly convex Banach spaces.
If $\delta _{B^{*}}(\epsilon )$ is modulus of the  convexity  of
space $_{B^{*}}$  and
$g_{B^{*}}(\epsilon ) = \delta _{B^{*}}(\epsilon )/\epsilon,
\; g^{-1}_{B^{*}}(\cdot )$
is an inverse function, then the following estimate
$$||Jx - Jy||_{B^{*}}
\leq C_1(g^{-1}_{B^{*}}{(2C_1L||x - y||)})$$.
is valid. Here $L$  is the constant from Lemma \ref{l1}.
\end {lem}
{\bf Proof}.
Analogously to (\ref {a15}), we can obtain the inequality
\begin{equation} \label{a7}
<Jx - Jy, x - y> \ge (2L)^{-1} \delta_{B^{*}} (||Jx-Jy||_{B^{*}}/C_2).
\end{equation}
for uniformly smooth space $B.$  From (\ref{a7}) one has
$$||Jx-Jy||_{B^{*}}||x - y|| \ge (2L)^{-1} \delta_{B^{*}}
(||Jx-Jy||_{B^{*}}/C_2)$$
Since $g_{B^{*}} (\epsilon)  = \delta_{B^{*}} (\epsilon)  / \epsilon $,
we can write
\begin{equation} \label{a8}
g_{B^{*}} (||Jx-Jy||_{B^{*}}/C_2) \leq 2C_2L||x - y||.
\end{equation}
This proves Lemma.

\section{Main results}
\setcounter{equation}{0}

{\bf The first problem}.  Estimate  $ ||P_{\Omega }x-P_{\Omega }y||$
via  $ ||x-y||$ .
\begin {th} \label{th1}
Let $B$ be the uniformly  convex  and uniformly smooth Banach spaces.
If $\; \delta _{B}(\epsilon )$ is modulus of the  convexity  of  space $B$
and
$g_{B}(\epsilon ) = \delta _{B}(\epsilon )/\epsilon,
g^{-1}_{B}(\cdot )$ is an inverse function, then
\begin{equation} \label{f1}
||P_{\Omega } x-P_{\Omega } y|| \leq  C g^{-1}_B  (2LC^2 g^{-1}_{B^{*}}
{(2CL||x - y||)})
\end{equation}
where $L $  is the constant from  Lemma \ref{l1}, and
$$ C = 2\max \lbrace 1,\; ||x-P_{\Omega }y||,\; ||y-P_{\Omega }x|| \rbrace.$$
\end{th}
{\bf Proof}. Let us denote $\bar x = P_{\Omega }x$ and $\bar y = P_{\Omega }y
$.
It is  known from (\ref{f0})  that
$$ <J(x - \bar x), \bar x - \bar y>  \; \ge \;  0,\;\;
<J(y - \bar y), \bar y - \bar x> \; \ge \; 0. $$
Therefore
\begin{equation} \label{f2}
<J(x - \bar y), \bar x - \bar y> \; \ge \;
   <J(x - \bar x) - J(x - \bar y), \bar y - \bar x>.
\end{equation}
Using Lemma \ref{l1} we obtain
\begin{equation} \label{f3}
<J(x - \bar x) - J(x - \bar y), \bar y - \bar x> \; \ge \;
(2L)^{-1}\delta _{B} (||\bar x - \bar y||/C_{1}).
\end{equation}
Here
$$ C_1 = 2 \max \lbrace  1, \sqrt{(||x - \bar x||^2
+ ||x - \bar y||^2)/2}\rbrace .$$
Now
$$ <J(x - \bar y) - J(y - \bar y), \bar x - \bar y> \; \ge \;
(2L)^{-1}\delta _{B} (||\bar x - \bar y||/C_{1}) $$
follows from (\ref {f2}) and (\ref {f3}).
Applying the Caushy-Schwarz inequality we get
\begin{equation} \label{f4}
g_{B}(||\bar x - \bar y||/C_{1}) \leq  2LC_{1}||J(x - \bar y)
- J(y - \bar y)||_{B^{*}}.
\end{equation}
In Lemma \ref{l2} the estimate
\begin{equation} \label{f5}
||J(x - \bar y) - J(y - \bar y)||_{B^{*}}
\leq C_2(g^{-1}_{B^{*}}{(2C_2L||x - y||)}).
\end{equation}
was obtained, where
$$ C_2 = 2 \max \lbrace  1, \sqrt {(||x - \bar y||^2
+ ||y - \bar y||^2)/2} \rbrace \leq  $$
$$2 \max \lbrace  1, ||x - \bar y||, ||y - \bar y|| \rbrace  \leq
2 \max \lbrace  1, ||x - \bar y||, ||y - \bar x|| \rbrace .$$
Thus, the inequality
$$||P_{\Omega } x - P_{\Omega } y|| \leq C_1 g^{-1}_{B}
(2LC_1C_2 g^{-1}_{B^{*}} (2C_2L||x - y||)) $$
is  realized  from (\ref {f4})   and  (\ref {f5}). Finally, we get
(\ref {f1})
because
$C_{1} \leq 2\max \lbrace  1, ||x - \bar y|| \rbrace .$

\begin{re}  It follows from (\ref {f1})  that  the  projection
operator $P_{\Omega }$ is
uniformly continuous on every bounded set of Banach space $B$.
\end{re}

\begin{re}  If $B$ and $B^*$ are  Hilbert  spaces $H$,  then the proof of
Theorem \ref{th1} gives
$||P_{\Omega } x-P_{\Omega } y||_H \leq ||x - y||_H .$ It is useful to remind
that $ \epsilon ^2 /8  \leq \delta _{H}(\epsilon ) \leq \epsilon ^2 /4.$
\end{re}

\begin{re}  If $y \in  \Omega $ then in (\ref {f1}) $C = 2 \max
\lbrace 1, 2 ||x-P_{\Omega }y|| \rbrace ;$ if $x \in  \Omega $,
then $C = 2 \max  \lbrace 1, 2 ||y-P_{\Omega }x|| \rbrace ,$ i.e.
$C = 2\max \lbrace 1, 2 ||x-y|| \rbrace .$
\end{re}

\begin{re}  Instead of (\ref {f5}) one  can  use  the  estimates  of  duality
mapping with gauge function too \cite{n,r}.
\end{re}

Let  $\Omega_1$ and  $\Omega_2$ be convex closed sets, $x \in B$ and
 $ H (\Omega_1, \Omega_2) \leq \sigma ,$ where
$$ H(\Omega_1, \Omega_2) =
\max \lbrace \sup_{z_1\in \Omega_1} \inf _{z_2\in \Omega_2}
 ||z_1 - z_2||, {\quad}\sup_{z_1\in \Omega_2} \inf _{z_2\in
\Omega_1}||z_1 - z_2|| \rbrace $$
is a Hausdorff distance between $\Omega_1$ and $\Omega_2.$ \\
\\
{\bf The second problem.}  Estimate  $||P_{\Omega _{1}}x-P_{\Omega _{2}}x||$
via $\sigma $.

\begin {th} \label{th2}
If $B$ is a uniformly convex space, $\delta _{B}(\epsilon )$ is  modulus  of
the convexity, and $\delta ^{-1}_{B}(\cdot )$ is an inverse function, then
\begin{equation} \label{f6}
||P_{\Omega _{1}}x-P_{\Omega _{2}}x||_{B} \leq
C_{1}\delta ^{-1}_{B}(4L(d+r) \sigma ),
\end{equation}
where $L\;$ is the constant $\;$ from Lemma \ref{l1},
$r = ||x||, \;d = \max \lbrace  d_{1},d_{2} \rbrace,\; d_{i} = $ {\it dist \/}
$(\theta, \Omega _{i}),\; i= 1,2,\; \theta $ is an origin of space $B$,
$C_{1} = 2\max \lbrace 1, r+d \rbrace .$
\end{th}
{\bf Proof}. Denote  $\bar x_1 = P_{\Omega _{1}}x,\; \bar x_2 =
P_{\Omega _{2}}x .$ From Lemma \ref {l1} we have
\begin{equation} \label{f7}
<J(x - \bar x_1) - J(x - \bar x_2), \bar x_2 - \bar x_1> \; \ge \;
(2L)^{-1}\delta _{B} (||\bar x_1 - \bar x_2||/C),
\end{equation}
$$C =  2 \max \lbrace  1, ||x - \bar x_1||, ||x - \bar x_2|| \rbrace .$$
There exist $\xi _{1} \in  \Omega _{1}$ such that
$||\bar x_{2}-\xi _{1}|| \leq \sigma $ and
$$ <J(x - \bar x_1), \bar x_2 - \bar x_1> =
   <J(x - \bar x_1), \bar x_2 - \xi _1> +
   <J(x - \bar x_1), \xi _1 - \bar x_1>  \leq \sigma ||x - \bar x_1|| ,$$
because $H (\Omega_1, \Omega_2) \leq \sigma ,$  and
$<J(x - \bar x_1), \xi _1 - \bar x_1>  \leq 0$ .

In the same way there exist $\xi _{2} \in  \Omega _{2}$ such  that
$||\bar x_{1}-\xi _{2}|| \leq \sigma $ and
$$ <J(x - \bar x_2), \bar x_1 - \bar x_2> =
   <J(x - \bar x_2), \bar x_1 - \xi _2> + $$
$$ <J(x - \bar x_2), \xi _2 - \bar x_2>  \leq \sigma ||x - \bar x_2|| ,$$
because $H (\Omega_1, \Omega_2) \leq \sigma ,$  and
$<J(x - \bar x_2), \xi _2 - \bar x_2>  \leq 0$.
Therefore
\begin{equation} \label{f8}
<J(x - \bar x_1) - J(x - \bar x_2), \bar x_2 - \bar x_1> \leq
\sigma (||x - \bar x_1|| + ||x - \bar x_2||)
\end{equation}
holds.  It is obvious that
$$||x - \bar x_1|| \leq  ||x - P_{\Omega _{1}} \theta||
\leq  ||x|| + || P_{\Omega _{1}} \theta|| \leq r + d ,$$
$$||x - \bar x_2|| \leq  ||x - P_{\Omega _{2}} \theta||
\leq  ||x|| + || P_{\Omega _{2}} \theta|| \leq r + d .$$
{}From (\ref {f7}) and (\ref {f8}) we can write
$$(2L)^{-1}\delta _{B} (||\bar x_1 - \bar x_2||/C_1) \leq 2 \sigma (r + d),\;
C_{1} = 2\max \;\lbrace 1, r+d \rbrace .$$
Hence the estimate (\ref {f6}) is valid.
\begin{re}  The more exact estimate (\ref {f6}) is
$$||P_{\Omega _{1}}x-P_{\Omega _{2}}x||_{B} \leq
C_{1}\delta ^{-1}_{B}(4LC_2 \sigma ), $$
$$C_1 =  2 \max \;\lbrace  1, ||x - \bar x_1||, ||x - \bar x_2|| \rbrace ,\;\;
C_2 =  2 \max \;\lbrace   ||x - \bar x_1||, ||x - \bar x_2|| \rbrace .$$
\end{re}

Let us notice for Hilbert space, that the estimate
$$||P_{\Omega _{1}}x-P_{\Omega _{2}}x||_{H} \leq
\sqrt {4 \sigma (2r + d) +  \sigma^2} $$
has been established in \cite{dk} (see also \cite{ol}).
However, applying the proof of Theorem \ref{th2} for
this  special case, one can obtain inequality
\begin{equation} \label{f9}
||P_{\Omega _{1}}x-P_{\Omega _{2}}x||_{H} \leq
\sqrt {2 \sigma (r + d)}.
\end{equation}
Indeed,
$$||\bar x_1 - \bar x_2||^2_H =
((x - \bar x_1) - (x - \bar x_2), \bar x_2 - \bar x_1)
\;\leq 2 \sigma (r + d) $$
takes place.
As one can see, the estimate (\ref {f9}) is obtained much more easily and of
a higher quality, than in \cite{dk} and \cite{ol}.\\
\\
Let $\Omega _{1}$  and $\Omega _{2}$ be convex closed sets, $x \in  B, y \in
B, H (\Omega _{1},\Omega _{2}) \leq  \sigma $. \\
\\
{\bf The third problem.}  Estimate  $ ||P_{\Omega _{1}}x-P_{\Omega _{2}}y||$
via  $ ||x-y||$ and $\sigma $.

The result follows immediately from (\ref {f1}) and (\ref {f6}):
$$ ||P_{\Omega _{1}}x-P_{\Omega _{2}}y|| \leq C g^{-1}_B
(2LC^2 g^{-1}_{B^{*}} {(2CL||x - y||)})  +
C_{1}\delta ^{-1}_{B}(4L(d+r) \sigma ), $$
where $r = ||y||, C = 2\max \;\lbrace 1,\;||x-P_{\Omega _{1}}y||,\;
||y-P_{\Omega _{1}}x|| \rbrace,$  and $C_{1}$ and $d$ are determined
in Theorem \ref{th2}.

\end{document}